\documentclass[12pt]{amsart}
\usepackage{amssymb}
\input{diagrams.tex}
\diagramstyle[scriptlabels,height=8mm,width=8mm]

\def\pdfsyncstart{}
\def\pdfsyncstop{}

\def\bdi{\pdfsyncstop\begin{diagram}}
\def\edi{\end{diagram}\pdfsyncstart}

\theoremstyle{plain}

\newtheorem{thm}{Theorem}[section]
\newtheorem{cor}[thm]{Corollary}
\newtheorem{lem}[thm]{Lemma}
\newtheorem{prop}[thm]{Proposition}

\theoremstyle{definition}
\newtheorem{defi}[thm]{Definition}
\newtheorem{defis}[thm]{Definitions}
\newtheorem{conj}[thm]{Problem}
\newtheorem{conv}[thm]{Convention}
\newtheorem{nota}[thm]{Notation}
\newtheorem{rem}[thm]{Remark}
\newtheorem{rems}[thm]{Remarks}
\newtheorem{exa}[thm]{Example}
\newtheorem{exas}[thm]{Examples}
\newtheorem{sit}[thm]{}

\newcommand{\brem}{\begin{rem}}
\newcommand{\brems}{\begin{rems}}
\newcommand{\erem}{\end{rem}}
\newcommand{\erems}{\end{rems}}
\newcommand{\bexa}{\begin{exa}}
\newcommand{\bexas}{\begin{exas}}
\newcommand{\eexa}{\end{exa}}
\newcommand{\eexas}{\end{exas}}
\newcommand{\bdefi}{\begin{defi}}
\newcommand{\edefi}{\end{defi}}
\newcommand{\bdefis}{\begin{defis}}
\newcommand{\edefis}{\end{defis}}
\newcommand{\bcor}{\begin{cor}}
\newcommand{\ecor}{\end{cor}}
\newcommand{\blem}{\begin{lem}}
\newcommand{\elem}{\end{lem}}
\newcommand{\bconv}{\begin{conv}}
\newcommand{\econv}{\end{conv}}
\newcommand{\bconj}{\begin{conj}}
\newcommand{\econj}{\end{conj}}
\newcommand{\bprop}{\begin{prop}}
\newcommand{\eprop}{\end{prop}}
\newcommand{\bthm}{\begin{thm}}
\newcommand{\ethm}{\end{thm}}
\newcommand{\bnota}{\begin{nota}}
\newcommand{\enota}{\end{nota}}
\newcommand{\bsit}{\begin{sit}}
\newcommand{\esit}{\end{sit}}
\newcommand{\be}{\begin{equation}}
\newcommand{\ee}{\end{equation}}
\newcommand{\bproof}{\begin{proof}}
\newcommand{\eproof}{\end{proof}}
\def\ba{\begin{array}}
\def\ea{\end{array}}
\def\bea{\begin{eqnarray}}
\def\eea{\end{eqnarray}}

\def\bnum{\begin{enumerate}}
\def\enum{\end{enumerate}}
\def\ext{{\rm ext}}

\newcommand{\Div}{\operatorname{Div}}

\newcommand{\Aut}{{\operatorname{Aut}}}

\def\cO{{\mathcal O}}

\def\bV{{\bar V}}
\def\bX{{\bar X}}

\newcommand{\A}{{\mathbb A}}

\newcommand{\PP}{{\mathbb P}}

\newcommand{\C}{{\mathbb C}}

\newcommand{\T}{{\mathbb T}}

\newcommand{\nlin}{\unitlength1mm\begin{picture}(0,9.25)
                     \put(0,0.75){\line(0,1){8.5}}
                    \end{picture}}

\newcommand{\vlin}[1]{\hspace{0.75mm}\unitlength1mm\begin{picture}
(#1,0)
                     \put(0,0){\line(1,0){#1}}
                    \end{picture}\hspace{0.75mm}\rule[-3mm]{0mm}{4mm}}

\newcommand{\lin}{\vlin{8.5}}

\newcommand{\cou}[2]{\unitlength1mm\begin{picture}(0,8)
  \put(0,0){\circle{1.5}}
  \put(0,3){\makebox(0,5)[b]{$#1$}}
  \put(0,-7){\makebox(0,4)[t]{$#2$}}
    \end{picture}
    \rule[-7mm]{0mm}{7mm}}

\newcommand{\crl}[2]{\unitlength1mm\begin{picture}(0,8)
  \put(0,0){\circle{1.5}}
  \put(-5,0){\makebox(0,5)[b]{$#1$}}
 \put(5,0){\makebox(0,5)[b]{$#2$}}
    \end{picture}
    \rule[-7mm]{0mm}{7mm}}

\newcommand{\cshiftup}[2]{\unitlength1mm\begin{picture}(0,9.25)
                     \put(0,10){\crl{#1}{#2}}
                    \end{picture}}

\title{On the Danilov-Gizatullin Isomorphism Theorem}

\author{Hubert Flenner}
\address{Fakult\"at f\"ur Mathematik,
Ruhr Universit\"at Bochum, Geb.\ NA 2/72, Universit\"ats\-str.\
150, 44780 Bochum, Germany}
\email{Hubert.Flenner@rub.de}

\author{Shulim Kaliman}
\address{Department of Mathematics,
University of Miami, Coral Gables, FL  33124, U.S.A.}
\email{kaliman@math.miami.edu}

\author{Mikhail Zaidenberg}
\address{Universit\'e
Grenoble I, Institut Fourier, UMR 5582 CNRS-UJF, BP 74, 38402 St.\
Martin d'H\`eres c\'edex, France} \email{zaidenbe@ujf-grenoble.fr}

\thanks{
{\bf Acknowledgements:} This research was done during a visit of
the first and second authors at the Institut Fourier, Grenoble and
of all three authors at the Max-Planck-Institute of Mathematics,
Bonn. They thank these institutions for the generous support and
excellent working conditions.}

\thanks{
\mbox{\hspace{11pt}}{\it 1991 Mathematics Subject Classification}:
14R05, 14R20.\\
\mbox{\hspace{11pt}}{\it Key words}: Hirzebruch surface, affine
surface}

\begin{document}

\begin{abstract}
A {\it Danilov-Gizatullin surface} is a normal affine surface
$V=\Sigma_d\setminus S$ which is a complement to an ample section
$S$ in a Hirzebruch surface $\Sigma_d$. By a surprising result due
to Danilov and Gizatullin [DaGi] $V$ depends  only on $n=S^2$ and
neither on $d$ nor on $S$. In this note we provide a new and
simple proof of this Isomorphism Theorem.
\end{abstract}

\maketitle

\section{The Danilov-Gizatullin theorem}

By definition, a Danilov-Gizatullin surface is the complement
$V=\Sigma_d\setminus S$ of an ample section $S$ in a Hirzebruch
surface $\Sigma_d$, $d\ge 0$.
In particular $n:=S^2>d$. The purpose of this note is to give
a short  proof of the following result of
Danilov and Gizatullin \cite[Theorem 5.8.1]{DaGi}.

\bthm\label{mthm} The isomorphism type of $V_n=\Sigma_d\setminus
S$ only depends on $n$. In particular, it neither depends on $d$
nor on the choice of the section $S$. \ethm

For other proofs we refer the reader to \cite{DaGi} and
\cite[Corollary 4.8]{CNR}. In the forthcoming paper 
\cite[Theorem 1.0.5]{FKZ3} we
extend the Isomorphism Theorem \ref{mthm} to a  larger
class of affine  surfaces. However, the proof of this latter
result is much harder.

\section{Proof of the Danilov-Gizatullin theorem}

\subsection{Extended divisors of Danilov-Gizatullin
surfaces}\label{2.1.0} Let as before $V=\Sigma_d\setminus S$ be a
Danilov-Gizatullin surface, where $S$ is  an ample section in a
Hirzebruch surface $\Sigma_d$, $d\ge 0$ with $n:=S^2>d$. Picking a
point, say, $A\in S$ and performing a sequence of $n$ blowups at
$A$ and its infinitesimally near points on $S$ leads to a new SNC
completion $(\bV, D)$ of $V$. The new boundary
$D=C_0+C_1+\ldots+C_n$ forms a {\em zigzag} i.e., a linear chain
of rational curves with weights $C_0^2=0$, $C_1^2=-1$ and
$C_i^2=-2$ for $i=2,\ldots, n$. Here $C_0\cong S$ is the proper
transform of $S$. The linear system $|C_0|$ on $\bV$ defines a
$\PP^1$-fibration $\Phi_0:\bV\to\PP^1$ for which $C_0$ is a fiber
and $C_1$ is a section. Choosing an appropriate affine coordinate on
$\PP^1=\A^1\cup \{\infty\}$ we may suppose that
$\Phi_0^{-1}(\infty)=C_0$ and $\Phi_0^{-1}(0)$ contains the
subchain $C_2+\ldots+C_n$ of $D$. The reduced curve
$D_\ext=\Phi_0^{-1}(0)\cup C_0\cup C_1$ is called the {\em
extended divisor} of the completion $(\bV,D)$ of $V$. The
following lemma appeared implicitly in the proof of Proposition 1
in \cite{Gi} (cf.\ also \cite{FKZ1}). To make this note
self-contained we provide a short argument.

\blem\label{exdi} \label{2.1} (a) For every $a\ne 0$ the fiber
$\Phi_0^{-1}(a)$ is reduced and isomorphic to $\PP^1$.

(b) $D_\ext=\Phi_0^{-1}(0)\cup C_0\cup C_1$ is an SNC divisor with
dual graph

\medskip

\be\label{graph}\;D_\ext:\qquad \cou{0}{C_0}\lin
\cou{-1}{C_1}\lin\cou{-2}{C_2}\lin \ldots \lin
\cou{\!\!\!\!\!\!\!\!\!\!\!\! -2}{C_s} \nlin\cshiftup{\!\!\!\!
1-s}{ F_1} \lin \ldots \lin \cou{\!\!\!\!\!\!\!\!\!\!\!\!\!
-2}{C_{n}}\nlin \cshiftup{-1}{F_0}\ee for some $s$ with $2\le s\le
n$. \elem

\bproof (a) follows easily from the fact that the affine surface
$V=\bV\backslash D$ does not contain complete curves.

To deduce (b), we note first that $\bV$ has Picard number $n+2$,
since $\bV$ is obtained from $\Sigma_d$ by a sequence of $n$
blowups. Since $C_1\cdot C_2=1$, the part $\Phi_0^{-1}(0)-C_2$ of
the fiber $\Phi_0^{-1}(0)$ can be blown down to a smooth point.
Since $C_1^2=-1$, after this contraction we arrive at the
Hirzebruch surface $\Sigma_1$, which has Picard number $2$. Hence
the fiber $\Phi_0^{-1}(0)$ consists of $n+1$ components. In other
words, $\Phi_0^{-1}(0)$ contains, besides the chain
$C_2+\ldots+C_n$, exactly 2 further components $F_0$ and $F_1$
called {\it feathers} \cite{FKZ1}. These are  disjoint smooth
rational curves, which meet the chain $C_2+\ldots+C_n$
transversally at two distinct smooth points. Indeed,
$\Phi_0^{-1}(0)$ is an SNC divisor without cycles and  the affine
surface $V$ does not contain complete curves. In particular,
$(F_0\cup F_1)\setminus D$ is a union of two disjoint smooth
curves on $V$ isomorphic to $\A^1$.

Since $\Phi_0^{-1}(0)-C_2$ can be blown down  to a smooth point
and $C_i^2=-2$  for $i\ge 2$, at least one of these feathers, call
it $F_0$, must be a $(-1)$-curve. We claim that $F_0$ cannot meet
a component $C_r$ with $3\le r\le n-1$. Indeed, otherwise the
contraction of $F_0+C_r+C_{r+1}$ would result in $C_{r- 1}^2=0$
without the total fiber over $0$ being irreducible, which is
impossible. Hence $F_0$ meets either $C_2$ or $C_n$.

If $F_0$ meets $C_2$ then $F_0+C_2+\ldots+C_n$ is contractible to
a smooth point. Thus the image of $F_1$ will become a smooth
fiber of the contracted surface. This is only possible if $F_1$ is
a $(-1)$-curve attached to $C_n$. Hence after interchanging $F_0$
and $F_1$ the divisor $D_\ext$ is as in (\ref{graph}) with $s=2$.

Therefore we may assume for the rest of the proof that $F_0$ is
attached at $C_n$ and $F_1$ at $C_s$, where $2\le s\le n$.
Contracting the chain $F_0+C_2+\ldots+C_n$ within the fiber
$\Phi_0^{-1}(0)$ yields an irreducible fiber $F'_1$ with
$(F'_1)^2=0$. This determines the index $s$ in a unique way,
namely, $s=1-F_1^2$. \eproof

\subsection{Jumping feathers} The construction in \ref{2.1.0}
depends on the initial choice of the point $A\in S$. In
particular, the extended divisor $D_\ext=D_\ext(A)$ and the
integer $s=s(A)$ depend on $A$. The aim of this subsection is to
show that $s(A)=2$ for a general choice of $A\in S$.

\bsit\label{2.2} Let $\bar F_0=\bar F_0(A)$ and $\bar F_1=\bar
F_1(A)$ denote the images of the feathers $F_0=F_0(A)$ and
$F_1=F_1(A)$, respectively, in the Hirzebruch surface $\Sigma_d$
under the blowdown $\sigma: \bV\to\Sigma_d$ of the chain
$C_1+\ldots+C_n$. These images meet each other and the original
section $S=\sigma(C_0)$ at the point $A$ and satisfy
\be\label{intersection} \bar F_0^2=0,\quad \bar F_0\cdot \bar
F_1=\bar F_0 \cdot S=1,\quad \bar F_1^2=n-2s+2, \quad \bar
F_1\cdot S=n-s+1\,, \ee where $s=s(A)$. Hence $\bar F_0=\bar
F_0(A)$ is the fiber through $A$ of the canonical projection
$\pi:\Sigma_d\to\PP^1$ and $\bar F_1=\bar F_1(A)$ is a section of
$\pi$. The sections $S$ and $\bar F_1$ meet only at $A$, where
they can be tangent (osculating).

We let below
 \be\label{lm} s_0=s(A_0)=\min_{A\in S}
\,\{s(A)\}\,,\quad l=\bar F_1(A_0)^2+1\,\,\,\mbox{and}\,\,\,
m=\bar F_1(A_0)\cdot S \,.\ee
For the next proposition see e.g.,
Lemma 7 and the following Remark in \cite{Gi}, or Proposition 4.8.11 in
\cite[II]{DaGi}. Our proof is based essentially on the same
idea.
\esit

\bprop\label{motherjump} (a) $s(A)=s_0$ for a general point $A\in
S$, and

(b) $s_0=2$. \eprop

\bproof For a general point $A\in
S$ and an arbitrary point $A'\in S$ we have $\bar F_1(A)\sim
\bar F_1(A')+k\bar F_0$ for some $k\ge 0$. Hence $\bar
F_1(A)^2=\bar F_1(A')^2+2k\ge \bar F_1(A')^2$. Using
(\ref{intersection}) it follows that
$$s(A)=1-F_1(A)^2\le s(A')=1-F_1(A')^2\,.$$ 
Thus $s(A)=s_0$ for all points $A$ in a Zariski open subset
$S_0\subseteq S$, which implies (a).

To deduce (b) we note that by (\ref{lm})
$$
l=n-2s_0+3\le n-s_0+1=m
$$
with equality if and only if $s_0=2$. Thus it is
enough to show that $l\ge m$. Restriction to $S$ yields
\be\label{m} \bar F_1(A)\vert S= m [A]\in\Div (S)\qquad\forall
A\in S_0\,. \ee Consider the linear systems
$$
|\bar F_1(A_0)|\cong \PP^l\quad\mbox{and}\quad |\cO_{S}
(m)|\cong \PP^m
$$
on $\Sigma_d$ and $S\cong \PP^1$,
respectively, and the linear map
$$\bdi\rho:\PP^l&\rDotsto &\PP^m\edi,\quad F\longmapsto F\vert S\,.$$
The set of divisors
$$\Gamma_m=\{m[A]\}_{A\in S}\,$$
represents a rational normal curve of degree $m$ in
$\PP^m=|\cO_{S} (m)|$. In view of (\ref{m}) the linear
subspace $\overline{\rho (\PP^l)}$ contains $\Gamma_m$.
Since the curve $\Gamma_m$ is linearly non-degenerate we have
$\overline{\rho (\PP^l)}=\PP^m$and so $l\ge m$, as desired.
\eproof

\subsection{Elementary shifts}
We consider as before a completion $V=\bV\setminus D$ of a
Danilov-Gizatullin surface $V$ as in \ref{2.1.0}.

\bsit\label{990}  Choosing $A$ generically, according to
Proposition \ref{motherjump} we may suppose in the sequel that
$s=s(A)=2$ and $F_0^2=F_1^2=-1$.

By (\ref{graph}) in Lemma \ref{2.1}, blowing down in $\bar V$ the
feathers $F_0,F_1$ and then the chain $C_3+\ldots+C_n$ yields the
Hirzebruch surface $\Sigma_1$, in which $C_0$ and $C_2$ become
fibers and $C_1$ a section. Reversing this contraction, the above
completion $\bV$ can be obtained  from $\Sigma_1$ by a sequence of
blowups as follows. The sequence starts by the blowup with center
at a point $P_3\in C_2\backslash C_1$ to create the next component
$C_3$ of the zigzag  $D$.
Then we perform subsequent blowups with centers at points $P_4,
\ldots, P_{n+1}$ infinitesimally near to $P_3$, where for each
$i=4,\ldots,n$ the blowup of $P_i\in C_{i-1}\backslash C_{i-2}$
creates the next component $C_i$ of the zigzag. The blowup with
center at $P_{n+1}\in C_n\backslash C_{n-1}$ creates the feather
$F_0$. Finally we blow up at a point $Q\in C_2\backslash C_1$
different from $P_3$ to create the feather $F_1$. In this way we
recover the given completion $\bV$ with extended divisor $D_\ext$
as in (\ref{graph}), where $s=2$.

We observe that the sequence  $P_3,\ldots, P_{n+1},
Q$ depends on the original triplet $(\Sigma_d,S,A)$. It follows
that, varying the points $P_3,\ldots, P_{n+1}, Q$
and then contracting the chain $C_1+\ldots+C_n=D-C_0$
on the resulting surface $\bar V$, we can obtain all possible
Danilov-Gizatullin surfaces
$$V=\bar V\setminus D\cong\Sigma_d \setminus
S\quad\mbox{with}\quad S^2=n\,\,\,\mbox{and}\,\,\, 0\le d\le
n-1\,.$$  Thus to deduce the Danilov-Gizatullin Isomorphism
Theorem \ref{mthm} it suffices to establish the following
fact.\esit

\bprop\label{2.6} The isomorphism type of the affine surface
$V=\bV\backslash D$ does not depend on the choice of the blowup
centers $P_3, \ldots, P_{n+1}$ and $Q$ as above. \eprop

The proof proceeds in several steps.

\bsit\label{qaff} First we note that in our construction it
suffices to keep track only of some partial completions rather
than of the whole complete surfaces. We  can choose affine
coordinates $(x,y)$ in $\Sigma_1\backslash (C_0\cup C_1)\cong
\A^2$ so that $C_2\backslash C_1=\{x=0\}$, $P=P_3=(0,0)$ and
$Q=(0,1)$. The affine surface $V$ can be obtained from the affine
plane $\A^2$ by performing subsequent blowups with centers at the
points
 $P_3,\ldots, P_{n+1}$ and $Q$ as in \ref{990} and then deleting
the curve $C_2\cup \ldots \cup C_n=D\backslash (C_0\cup C_1)$.

With $X_2=\A^2$, for every $i= 3,\ldots, n+1$ we let $X_i$ denote
the result of the subsequent blowups of $\A^2$ with centers
$P_3, \ldots, P_i$. This gives a tower of blowups
\be\label{tower}
\bar V\backslash (C_0\cup C_1)=:X_{n+2}\to
X_{n+1}\to X_n\to\ldots\to X_2=\A^2\,,
\ee
where in the last step the point
$Q$ is blown up to create $F_1$.
\esit

\bsit\label{aut} Let us exhibit a special case of this
construction. Consider the standard action
$$(\lambda_1,\lambda_2):(x,y)\mapsto (\lambda_1 x,\lambda_2 y)$$
of the 2-torus $\T=(\C^*)^2$ on the affine plane $X_2=\A^2$. We
claim that there is a unique sequence of points $(0,0)=P_3=P_3^o,
\ldots ,P_{n+1}=P^o_{n+1}$ such that the torus action can be
lifted to $X_i$ for $i=3,\ldots ,n+1$. Indeed, if by induction the
$\T$-action is lifted already to $X_{i}$ with $i\ge 2$, then on
$C_i\backslash C_{i-1}\cong \A^1$ the induced $\T$-action has a
unique fixed point $P_{i+1}^o$. Blowing up this point the
$\T$-action can be lifted further to $X_{i+1}$. Blowing up finally
$Q=(0,1)\in C_2\setminus C_1$ and deleting $C_2\cup \ldots\cup
C_n$ we arrive at a unique {\em standard} Danilov-Gizatullin
surface $V_{\rm st}=V_{\rm st}(n)$.

Let us note that $\T$ acts transitively on 
$(C_2\setminus C_1)\setminus \{(0,0)\}$. Thus up to
isomorphism, 
the resulting affine surface $V_{\rm st}$ 
does not depend on the choice of $Q$.
\esit

\bsit\label{aut1} Consider now an automorphism $h$ of $\A^2$
fixing the $y$-axis pointwise. It moves the blowup centers $P_4,
\ldots, P_{n+1}$ to new positions $P_4', \ldots, P_{n+1}'$, while
$P_3$ and $Q$ remain unchanged. It is easily seen that $h$ induces
an isomorphism between $V$ and the resulting new affine surface
$V'$. We show in Lemma \ref{2.7} below that applying  a suitable
automorphism $h$, we can choose $V'$ to be the standard surface
$V_{\rm st}$ as in \ref{aut}. This implies immediately Proposition
\ref{2.6} and as well Theorem \ref{mthm}. More precisely, our $h$
will be composed of {\em elementary shifts} \be\label{elmap}
h_{a,t}: (x,y)\mapsto (x,y+ax^{t}),\quad\mbox{where}\quad a\in
\C\quad\mbox{and}\quad t\ge 0\,. \ee \esit

\blem\label{2.7} By a sequence of elementary shifts as in
(\ref{elmap}) we can move the blowup centers $P_4, \ldots ,P_n$
into the points $P_4^o, \ldots, P_n^o$ so that $V$ is isomorphic
to $V_{\rm st}$. \elem

\bproof Since $X_2=\A^2$ the assertion is obviously true for
$i=2$. The point $P_3=(0,0)$ being fixed by $\T$, the torus
action can be lifted to $X_3$. The blowup with center at
$P_3$ has a coordinate presentation
$$
(x_3, y_3)=(x, y/x)\,,\quad\mbox{or, equivalently,}\quad
(x,y)=(x_3, x_3y_3)\,,
$$
where the exceptional curve $C_3$ is given by $x_3=0$ and the
proper transform of $C_2$ by $y_3=\infty$. The action of $\T$ in
these coordinates is
$$
(\lambda_1,\lambda_2).(x_3,y_3)=(\lambda_1x_3,\lambda_1^{-1}\lambda_2y_3)\,,
$$
while the elementary shift $h_{a,t}$ can be written as
\be\label{transl} h_{a,t}: (x_3,y_3)\mapsto (x_3,
y_3+ax_3^{t-1})\,. \ee Thus in  $(x_3,y_3)$-coordinates
$P_4^o=(0,0)$. Furthermore for $t=1$, the shift $h_{a,1}$ yields a
translation on the axis $C_3\backslash C_2=\{x_3=0\}$, while
$h_{a,t}$ with $t\ge 2$ is the identity on this axis. Applying $h_{a,1}$
for a suitable $a$ we can move the point $P_4\in C_3\backslash
C_2$ to $P_4^o$. Repeating the argument recursively, we can
achieve that $P_{i}=P_{i}^o$ for $i\le n+1$, as required. \eproof

\brems
 1. The surface $X_{n+1}$ as in \ref{aut} is toric,
and the $\T$-action on $X_{n+1}$ stabilizes  the chain
$C_2\cup\ldots\cup C_n\cup F_0$. There is a 1-parameter subgroup
$G$ of the torus (namely, the stationary subgroup of the point
$Q=(0,1)$), which lifts to $X_{n+2}$ and then restricts to $V_{\rm
st}=X_{n+2}\setminus (C_2\cup\ldots\cup C_n)$. Fixing an isomorphism $G\cong
\C^*$ gives a $\C^*$-action on $V_{\rm st}$. As follows from
\cite[1.0.6]{FKZ3}, $V_{\rm st}=V_{\rm st}(n)$ is the normalization of the
surface $W_n\subseteq\A^3$ with equation
$$
x^{n-1}y=(z-1)(z+1)^{n-1}\,.
$$
For $n\ge 3$ this surface has
non-isolated singularities, and is equipped with the
$\C^*$-action $\lambda.(x,y,z)=(\lambda x, \lambda^{n-1}y, z)$.
Due to the Danilov-Gizatullin Isomorphism Theorem \ref{mthm},
any Danilov-Gizatullin surface $V_n$ is isomorphic to the
normalization of $W_n$.

2. However, the specific $\C^*$-action on $V_n$ obtained in this
way is not unique as was observed by Peter Russell. According
to Proposition 5.14 in \cite{FKZ1},  in $\Aut (V_n)$ there are
exactly $n-1$ different conjugacy classes of such actions
corresponding to different choices of $s=2,\ldots,n$ in 
diagram (\ref{graph}).
Let us sketch a construction of these classes which does not rely
on DPD-presentations as in {\em loc.cit}, but follows a procedure
similar to those used in the proof above.

Given $s\in \{2,\ldots,n\}$, starting with
$\bX_2=\Sigma_1\to\PP^1$ and a chain $C_0+C_1+C_2$ on $\Sigma_1$
as in \ref{990} and \ref{qaff}, we blow up the point $(0,0)\in
C_2$ creating the feather $F_1$, then at the point $C_2\cap F_1$
creating $C_3$ etc., until the component $C_s$ is created. The
standard torus action on $\Sigma_1$ lifts to the resulting surface
$\bX_{s+1}$ stabilizing the linear chain $F_1+C_0+\ldots+C_s$.
Next we blowup at a point $P\in C_s\setminus (F_1\cup C_{s-1})$
creating a new component $C_{s+1}$, and we lift the action of the
1-parameter subgroup $G=$Stab$_P(\T)$ to the resulting surface
$\bX_{s+2}$. Choosing an appropriate isomorphism $G\cong\C^*$ we
may assume that $C_s$ is attractive for the resulting
$\C^*$-action $\Lambda_s$ on $\bX_{s+2}$.  We continue blowing up
subsequently at the fixed points of this action on the curves
$C_{i+1}\setminus C_i$, $i=s,\ldots,n$ creating components
$C_{s+2}, \ldots,C_n$ and the feather $F_0$. Finally we arrive at
a $\C^*$-surface $\bV=\bX_{n+2}$ with an extended divisor as in
(\ref{graph}). Contracting $C_1+\ldots +C_n$ exhibits the open
part $V=\bV\setminus D$, where $D=C_0+\ldots +C_n$, as a
complement to an ample section in a Hirzebruch surface. Thus
$V=V_n$ is a Danilov-Gizatullin surface of index $n$ endowed
with a $\C^*$-action say, $\Lambda_s$, such that $\bV$ is its
equivariant standard completion. Note that the isomorphism class
of $(\bV, D)$ is independent on the choice of the point $P\in C_s\setminus
(F_1\cup C_{s-1})$. Indeed this point can be moved by the
$\T$-action yielding conjugated $\C^*$-actions on $V_n$.

Contracting the chain $C_1+\ldots +C_n$ leads to a Hirzebruch
surface $\Sigma_d$ such that the image of $F_0$ is a fiber of the ruling
$\Sigma_d\to\PP^1$. Moreover, the image $S$ of $C_0$ is an ample
section with $S^2=n$ so that $V_n=\Sigma_d\setminus S$.
The image of $F_1$ is another section with $F_1^2=n+2-2s$.
In particular, if this number is negative then $d= 2s-2-n$.

One can show that the $\Lambda_s$, $s=2,\ldots,n$ represent all conjugacy
classes
of $\C^*$-actions on $V_n$. Moreover, inverting
the action $\Lambda_s$ with respect to the isomorphism $t\mapsto
t^{-1}$ of $\C^*$ yields the action $\Lambda_{n-s+2}$. Thus after
inversion, if necessary, we may suppose that $2s-2\ge n$ 
so that $V_n\cong \Sigma_d\backslash S$ as above with $d=2s-2-n$.

%One can show that the $\Lambda_s$, $s=2,\ldots,n$ are up to
%conjugation all different $\C^*$-actions on $V_n$. Moreover, inverting
%the action $\Lambda_s$ with respect to the isomorphism $t\mapsto
%t^{-1}$ of $\C^*$ yields the action $\Lambda_{n-s+2}$. Thus after
%inversion, if necessary, we may suppose that $2s-2\ge n$. 
%Contracting the chain $C_1+\ldots +C_n$ leads to a Hirzebruch
%surface $\Sigma_d\to\PP^1$, where $d=2s-2-n$, with an ample
%section $S\cong C_0$ ($S^2=n$), the exceptional
%section $F_1$ and a fiber $F_0$ through a point $A\in F_1\cap S$,
%so that $V_n=\Sigma_d\setminus S$.

%It is easy to check that the feathers $F_0,\,F_1$ are fiber
%components of the orbit map $V_n\to\A^1$, where $F_1$ is of
%multiplicity 1, while $F_0$ has multiplicity $n-s+1$ in the fiber.
%Thus $s$ is an invariant of the conjugacy class of the resulting
%$\C^*$-action on $V_n$. We obtain in this way  $n-1$ different
%such conjugacy classes corresponding to $s=2,\ldots,n$.

%Conversely, given a $\C^*$-action on $V_n$ there exists an
%equivariant standard completion $(\bV,D)$ of $V_n$ with extended
%divisor as in (\ref{graph}) for a certain $s\in \{2,\ldots,n\}$,
%see \cite{FKZ1}. Reversing the above construction it is easily
%seen that $\bV$ with the given $\C^*$-action can be actually
%obtained via the above process. Hence this $\C^*$-action on $V_n$
%belongs to the corresponding conjugacy class with number $s$.

3. As was remarked by Peter Russell, with the exception of
Proposition \ref{motherjump} our proof is also valid for
Danilov-Gizatullin surfaces over an algebraically closed field of
any characteristic $p$. Moreover Proposition \ref{motherjump}
holds as soon as $p=0$ or $p$ and $m$ are coprime. In particular
it follows that the Isomorphism Theorem holds in the cases $p=0$
and $p\ge n-2$. This latter
% also this more general
result was shown already in \cite{DaGi}. However for $p=2$ and
$n=56$ there is an infinite number of isomorphism types of 
Danilov-Gizatullin surfaces; see \cite[\S 9]{DaGi}. \erems

%Indeed by Theorem 2.9(b) in
%\cite{FKZ1}, given a $\C^*$-action on $V_n$ there is an
%equivariant completion of $V_n$ as in (\ref{graph}) with a feather
%$F_1$ attached at the component $C_s$ for some $s\in [2,n]$.
%Blowing down $C_1,F_0,C_n,\ldots,C_3, F_1$ is equivariant and
%yields the projective plane $\PP^2$ with a linear $\C^*$-action
%which stabilizes two lines, images of $C_0$ and $C_2$.

%Reversing this procedure, first we blow up subsequently at fixed
%points of the torus action to create the curves $C_1$ and
%$F_1+C_3+\ldots+C_s$. There is a $1$-parameter subgroup of the
%torus $\T$ which fixes $C_s$ pointwise. Next we choose an
%arbitrary point $p\in C_s\setminus C_{s-1}$ to create $C_{s+1}$,
%and the other blowups are done at unique fixed points of the
%latter $\C^*$-action. The choice of the point $p$ is irrelevant
%since, at step $s$, the toric action is non-trivial on the
%component $C_s$, hence it moves the center $p$ so that the final
%$\C^*$-surfaces are all equivariantly isomorphic. \cs

\end{document}